\documentclass{article}
\usepackage{graphicx} 
\usepackage[dvipsnames]{xcolor}

\usepackage{tikz}
\usepackage{circuitikz}
\usepackage{tikz-3dplot}
\usetikzlibrary{arrows.meta}
\usetikzlibrary{positioning}

\usepackage{amsbsy}
\usepackage{graphicx}
\usepackage{booktabs}
\usepackage{tabularx}
\usepackage{lipsum}
\usepackage{amsthm} 
\usepackage{cite}
\usepackage{amsthm}
\usepackage{mathtools}
\usepackage{amssymb}
\usepackage{amsmath}
\usepackage{algorithm}
\usepackage{algorithmic}
\usepackage{authblk}

\newtheorem{fact}{Fact} 

\title{Extremal poker hand rankings: why the standard 52 card deck and a 3044 card deck are special}

\author{Christopher Williamson}

\affil{\texttt{cw.@williamsonchris.com}}

\date{September 2025}

\begin{document}

\maketitle

\abstract{

We study poker hand rankings in the partially generalised setting of a deck with $r$ ranks, rather than the typical 13 ranks. We provide the hand rankings for all $r$ and observe some interesting phenomena such as the smallest $r$ such that flushes rank below one-pair hands. Perhaps surprisingly, as $r$ grows without bound, the hand ranking is not stable until $r=761$ (a 3044 card deck). We consider showdown frequency, which is the frequency that a given type of hand is declared by a player at showdown, and make note of counterintuitive instances in which a hand with lower absolute frequency than some other hand nonetheless has a higher showdown frequency. This can be interpreted as a form of Gadbois paradox but in the typical setting of poker without wild cards. Conveniently, the standard deck with 13 ranks turns out to be the smallest deck that avoids a discrepancy between absolute frequency and showdown frequency for all hand types other than having a high card. 

}

\section{Introduction}
Texas hold 'em poker is typically played with a standard deck of 52 playing cards, comprising 13 ranks (2 through Ace) in each of four suits (spades, hearts, clubs, diamonds). Players form five-card hands from a set of seven cards (two personal cards and five cards on the board). These hands are then ranked in an order determined by their rarity. Whereas the standard ranking of poker hands is well-established, it is calibrated specifically for a standard 52-card deck. When the deck composition changes, the frequencies of different hand types also change, along with the relative rankings. 

Short-deck poker is an increasingly popular poker variant that uses a deck where the ranks 2-5 have been removed. As a result, flushes rank higher than full houses. We compute frequency-based hand rankings as the number of card ranks $r$ grows without bound. Each deck has the standard suits and thus $4r$ total cards.

Our calculations also have a relationship with the Gadbois paradox. In 1996, Gadbois~\cite{gadbois1996poker} (see also~\cite{emert1996inconsistencies}) showed that the introduction of a wild card can create ``paradoxes'' in the task of creating a hand ranking. In a nutshell, Gadbois showed that if two-pair hands are ranked below three-of-a-kind (as is typical), then two-pair hands will become more rare than three-of-a-kind when one considers that, players will form a three-of-a-kind hand rather than a two-pair hand whenever they have a choice. Moreover, if one redefines the hand rankings to reflect that two-pair is rarer than three-of-a-kind, then two-pair becomes more common again as players will now use their wild cards to form two-pair hands over three-of-a-kind when they have that option. 

What Gadbois shows is that hand types that are rarer in an absolute sense may nonetheless be declared at showdown more often than more common hands. This stems from the player's choice of assigning a wild card value. But this phenomenon does not occur only in the presence of wild cards. Even without wild cards, players have a choice of which five of seven cards to use to form a hand and have a choice (which will be motivated by a previously-defined hand ranking) of what hand to declare at show-down. One could imagine that in a deck with very few ranks, conditioned on holding a hand with three-of-a-kind, that one is more likely than not to have the option to declare a full house once all cards have appeared. In this case, under the typical ranking, a full-house would be declared at showdown more often than three-of-a-kind, despite the fact that every full house contains a three-of-a-kind and thus is of a lower absolute frequency. We study this and show for which values of $r$ and which hand types such a wild-card-free Gadbois phenomenon exists.

\section{Preliminaries}
We assume knowledge of poker hand definitions and use the abbreviations HC, 1P, 2P, 3X, ST, FL, FH, 4X, SF to denote high-card, one-pair, two-pair, three-of-a-kind, straight, flush, full house, four-of-a-kind, and straight flush, respectively.\footnote{The only controversial choice we make is to consider 4X as containing a 2P. And naturally, a royal flush is merely a type of straight flush.}

Our first goal is to determine hand rankings for general values of $r$. We compute hand frequencies in an inclusive way that is agnostic to any particular hand ranking. For example, we cannot count a hand that contains a straight and a flush as a flush but not a straight, because this is to pre-suppose that flushes are the higher-ranking hand.\footnote{In fact, this inclusive counting method is proposed by~\cite{emert1996inconsistencies} to resolve the Gadbois paradox.} 

When we write $freq_r(h)$ for some $h\in\{\text{HC, 1P, 2P, 3X, ST, FL, FH, 4X, SF}\}$, we mean the number of hands, out of the total set of size $\binom{4r}{7}$ that contain the hand type $h$. Because we count inclusively, we will have that $\sum_{h}freq_r(h)>\binom{4r}{7}$. In mathematical notation, we let $S_r$ be the set of all seven card hands from a deck with $4r$ cards (four suits and $r$ ranks). For an element $s\in S_r$, we say that $h\in s$ if there exists a subset of five cards of the seven in $s$ such that those five cards contain a hand of type $h$. Then we have:
\[
freq_r(h)=|\{s\in S_r|h\in s\}|
\]

\section{Hand rankings}
We computed $freq_r$ for all hand types and varying values of $r$. With these frequencies in hand, we then derive the \emph{frequency ranking} where we rank the hand types so that type $h_0$ is ranked above $h_1$ exactly when $freq_r(h_0) < freq_r(h_1)$. Closed form expressions for $freq_r$ and their asymptotic behaviour are provided in Table 2; here we provide the final hand rankings. 

To reduce clutter in the diagram, we omit HC and SF hands, as these are always the lowest (resp., highest) ranking hand types (note that when $r=5$, every FL is also a SF). 

\begin{figure}[!ht]
\centering
\resizebox{1\textwidth}{!}{%
\begin{circuitikz}
\tikzstyle{every node}=[font=\footnotesize]
\draw [->, >=Stealth] (3.75,8.25) -- (23.5,8.25);
\draw  (5,12.75) circle (0.5cm) node {\normalsize 3X} ;
\draw  (5,14.25) circle (0.5cm) node {\normalsize FH} ;
\draw  (5,15.75) circle (0.5cm) node {\normalsize ST} ;
\draw  (5,17.25) circle (0.5cm) node {\normalsize 4X} ;
\draw  (5,18.75) circle (0.5cm) node {\normalsize FL} ;
\draw  (5,11.25) circle (0.5cm) node {\normalsize 2P} ;
\draw  (5, 9.75) circle (0.5cm) node {\normalsize 1P} ;

\draw  (7,17.25) circle (0.5cm) node {\normalsize 4X} ;
\draw  (7,18.75) circle (0.5cm) node {\normalsize FL} ;
\draw  (7,15.75) circle (0.5cm) node {\normalsize FH} ;
\draw  (7,14.25) circle (0.5cm) node {\normalsize 3X} ;
\draw  (7,12.75) circle (0.5cm) node {\normalsize ST} ;
\draw  (7,11.25) circle (0.5cm) node {\normalsize 2P} ;
\draw  (7, 9.75) circle (0.5cm) node {\normalsize 1P} ;

\draw  (9,18.75) circle (0.5cm) node {\normalsize 4X} ;
\draw  (9,17.25) circle (0.5cm) node {\normalsize FL} ;
\draw  (9,15.75) circle (0.5cm) node {\normalsize FH} ;
\draw  (9,14.25) circle (0.5cm) node {\normalsize 3X} ;
\draw  (9,12.75) circle (0.5cm) node {\normalsize ST} ;
\draw  (9,11.25) circle (0.5cm) node {\normalsize 2P} ;
\draw  (9, 9.75) circle (0.5cm) node {\normalsize 1P} ;

\draw  (11,12.75) circle (0.5cm) node {\normalsize 3X} ;
\draw  (11,15.75) circle (0.5cm) node {\normalsize FH} ;
\draw  (11,17.25) circle (0.5cm) node {\normalsize FL} ;
\draw  (11,18.75) circle (0.5cm) node {\normalsize 4X} ;
\draw  (11,14.25) circle (0.5cm) node {\normalsize ST} ;
\draw  (11,11.25) circle (0.5cm) node {\normalsize 2P} ;
\draw  (11, 9.75) circle (0.5cm) node {\normalsize 1P} ;

\draw  (13,12.75) circle (0.5cm) node {\normalsize 3X} ;
\draw  (13,14.25) circle (0.5cm) node {\normalsize ST} ;
\draw  (13,15.75) circle (0.5cm) node {\normalsize FL} ;
\draw  (13,17.25) circle (0.5cm) node {\normalsize FH} ;
\draw  (13,18.75) circle (0.5cm) node {\normalsize 4X} ;
\draw  (13,11.25) circle (0.5cm) node {\normalsize 2P} ;
\draw  (13, 9.75) circle (0.5cm) node {\normalsize 1P} ;

\draw  (15,18.75) circle (0.5cm) node {\normalsize 4X} ;
\draw  (15,17.25) circle (0.5cm) node {\normalsize FH} ;
\draw  (15,15.75) circle (0.5cm) node {\normalsize ST} ;
\draw  (15,14.25) circle (0.5cm) node {\normalsize FL} ;
\draw  (15,12.75) circle (0.5cm) node {\normalsize 3X} ;
\draw  (15,11.25) circle (0.5cm) node {\normalsize 2P} ;
\draw  (15, 9.75) circle (0.5cm) node {\normalsize 1P} ;

\draw  (17,18.75) circle (0.5cm) node {\normalsize 4X} ;
\draw  (17,17.25) circle (0.5cm) node {\normalsize FH} ;
\draw  (17,15.75) circle (0.5cm) node {\normalsize ST} ;
\draw  (17,14.25) circle (0.5cm) node {\normalsize 3X} ;
\draw  (17,12.75) circle (0.5cm) node {\normalsize FL} ;
\draw  (17,11.25) circle (0.5cm) node {\normalsize 2P} ;
\draw  (17, 9.75) circle (0.5cm) node {\normalsize 1P} ;

\draw  (19,18.75) circle (0.5cm) node {\normalsize 4X} ;
\draw  (19,17.25) circle (0.5cm) node {\normalsize ST} ;
\draw  (19,15.75) circle (0.5cm) node {\normalsize FH} ;
\draw  (19,14.25) circle (0.5cm) node {\normalsize 3X} ;
\draw  (19,12.75) circle (0.5cm) node {\normalsize 2P} ;
\draw  (19,11.25) circle (0.5cm) node {\normalsize FL} ;
\draw  (19, 9.75) circle (0.5cm) node {\normalsize 1P} ;

\draw  (21,18.75) circle (0.5cm) node {\normalsize 4X} ;
\draw  (21,17.25) circle (0.5cm) node {\normalsize ST} ;
\draw  (21,15.75) circle (0.5cm) node {\normalsize FH} ;
\draw  (21,14.25) circle (0.5cm) node {\normalsize 3X} ;
\draw  (21,12.75) circle (0.5cm) node {\normalsize 2P} ;
\draw  (21,11.25) circle (0.5cm) node {\normalsize 1P} ;
\draw  (21, 9.75) circle (0.5cm) node {\normalsize FL} ;

\draw  (23,18.75) circle (0.5cm) node {\normalsize ST} ;
\draw  (23,17.25) circle (0.5cm) node {\normalsize 4X} ;
\draw  (23,15.75) circle (0.5cm) node {\normalsize FH} ;
\draw  (23,14.25) circle (0.5cm) node {\normalsize 3X} ;
\draw  (23,12.75) circle (0.5cm) node {\normalsize 2P} ;
\draw  (23,11.25) circle (0.5cm) node {\normalsize 1P} ;
\draw  (23, 9.75) circle (0.5cm) node {\normalsize FL} ;

\draw [short] (5.5,18.75) -- (6.5,18.75);
\draw [short] (5.5,17.25) -- (6.5,17.25);
\draw [short] (5.5,15.75) -- (6.5,12.75);
\draw [short] (5.5,14.25) -- (6.5,15.75);
\draw [short] (5.5,12.75) -- (6.5,14.25);
\draw [short] (7.5,14.25) -- (8.5,14.25);
\draw [short] (7.5,12.75) -- (8.5,12.75);
\draw [short] (9.5,14.25) -- (10.5,12.75);
\draw [short] (9.5,12.75) -- (10.5,14.25);
\draw [short] (7.5,15.75) -- (8.5,15.75);
\draw [short] (7.5,17.25) -- (8.5,18.75);
\draw [short] (7.5,18.75) -- (8.5,17.25);
\draw [short] (9.5,18.75) -- (10.5,18.75);
\draw [short] (11.5,18.75) -- (12.5,18.75);
\draw [short] (13.5,18.75) -- (14.5,18.75);
\draw [short] (9.5,17.25) -- (10.5,17.25);
\draw [short] (9.5,15.75) -- (10.5,15.75);
\draw [short] (11.5,17.25) -- (12.5,15.75);
\draw [short] (11.5,15.75) -- (12.5,17.25);
\draw [short] (13.5,17.25) -- (14.5,17.25);
\draw [short] (11.5,14.25) -- (12.5,14.25);
\draw [short] (13.5,14.25) -- (14.5,15.75);
\draw [short] (13.5,15.75) -- (14.5,14.25);
\draw [short] (11.5,12.75) -- (12.5,12.75);
\draw [short] (13.5,12.75) -- (14.5,12.75);

\draw [short] (17.5,18.75) -- (18.5,18.75);
\draw [short] (19.5,18.75) -- (20.5,18.75);
\draw [short] (21.5,18.75) -- (22.5,17.25);
\draw [short] (21.5,17.25) -- (22.5,18.75);
\draw [short] (19.5,17.25) -- (20.5,17.25);
\draw [short] (17.5,15.75) -- (18.5,17.25);
\draw [short] (17.5,17.25) -- (18.5,15.75);
\draw [short] (19.5,15.75) -- (20.5,15.75);
\draw [short] (21.5,15.75) -- (22.5,15.75);

\draw [short] (21.5,14.25) -- (22.5,14.25);
\draw [short] (19.5,14.25) -- (20.5,14.25);
\draw [short] (17.5,14.25) -- (18.5,14.25);

\draw [short] (19.5,12.75) -- (20.5,12.75);
\draw [short] (21.5,12.75) -- (22.5,12.75);
\draw [short] (21.5,11.25) -- (22.5,11.25);
\draw [short] (21.5,9.75) -- (22.5,9.75);

\draw [short] (5.5,9.75) -- (6.5,9.75);
\draw [short] (7.5,9.75) -- (8.5,9.75);
\draw [short] (9.5,9.75) -- (10.5,9.75);
\draw [short] (11.5,9.75) -- (12.5,9.75);
\draw [short] (13.5,9.75) -- (14.5,9.75);
\draw [short] (15.5,9.75) -- (16.5,9.75);
\draw [short] (17.5,9.75) -- (18.5,9.75);

\draw [short] (19.5,9.75) -- (20.5,11.25);
\draw [short] (19.5,11.25) -- (20.5,9.75);
\draw [short] (17.5,12.75) -- (18.5,11.25);
\draw [short] (17.5,11.25) -- (18.5,12.75);

\draw [short] (5.5,11.25) -- (6.5,11.25);
\draw [short] (7.5,11.25) -- (8.5,11.25);
\draw [short] (9.5,11.25) -- (10.5,11.25);
\draw [short] (11.5,11.25) -- (12.5,11.25);
\draw [short] (13.5,11.25) -- (14.5,11.25);
\draw [short] (15.5,11.25) -- (16.5,11.25);

\node [font=\normalsize] at (5,7.75) {5};
\node [font=\normalsize] at (7,7.75) {6};
\node [font=\normalsize] at (9,7.75) {7};
\node [font=\normalsize] at (11,7.75) {8-12};
\node [font=\normalsize] at (13,7.75) {13,14};
\node [font=\normalsize] at (15,7.75) {15-18};
\node [font=\normalsize] at (17,7.75) {19-35};
\node [font=\normalsize] at (19,7.75) {36-306};
\node [font=\normalsize] at (21,7.75) {307-760};
\node [font=\normalsize] at (23,7.75) {761+};
\node [font=\normalsize, rotate around={90:(0,0)}] at (3,14) {Frequency ranking};
\node [font=\footnotesize] at (13,7.25) {Standard deck};
\node [font=\footnotesize] at (11,7.25) {Short deck};
\node [font=\footnotesize] at (16.5,10.75) {};
\draw [short] (15.5,12.75) -- (16.5,14.25);
\draw [short] (15.5,14.25) -- (16.5,12.75);
\draw [short] (15.5,15.75) -- (16.5,15.75);
\draw [short] (15.5,17.25) -- (16.5,17.25);
\draw [short] (15.5,18.75) -- (16.5,18.75);
\draw [short] (3.75,8.25) -- (3.75,20);
\node [font=\normalsize] at (14.5,6.25) {Number of ranks $r$};

\end{circuitikz}
}%

\label{fig:my_label}
\caption{Frequency rankings for general $r$. We do not consider $r<5$ as this makes some hand types impossible. The cases $r=9$ and $r=13$ correspond to short-deck poker and standard poker, respectively.}
\end{figure}
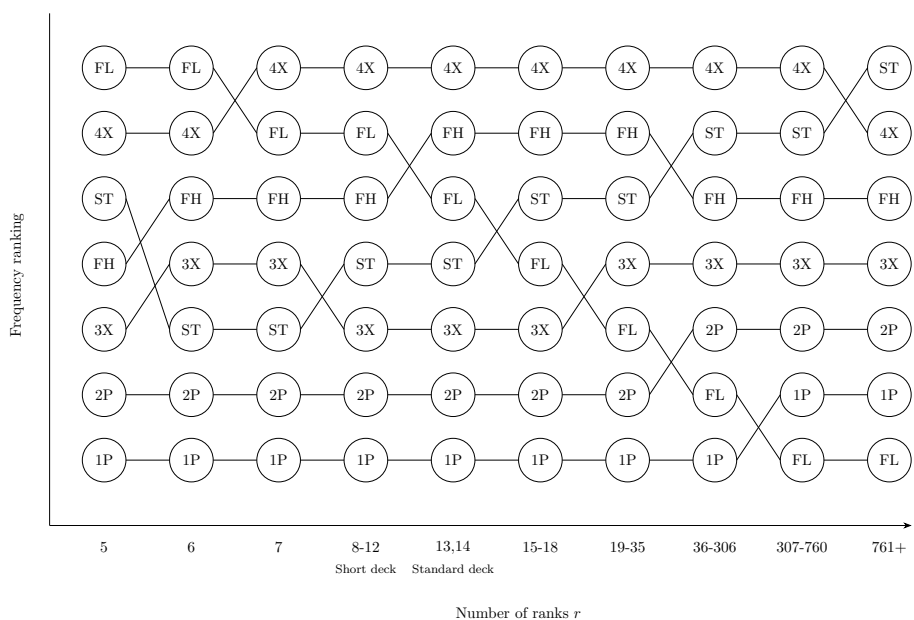

On the basis of these frequency rankings, we list a few facts:

\begin{fact}
Flushes eventually rank below one-pair hands, although this only happens when $r\geq 307$.
\end{fact}

\begin{fact}
Straights eventually become the rarest type of hand (other than straight-flushes), although this only happens when $r\geq 761$.
\end{fact}

\begin{fact}
The frequency ranking is stable for all $r\geq 761$ as $r\rightarrow \infty$.
\end{fact}

To prove Fact 3, one can expand out the formulas in Table 2 as polynomials and, for each pair of hand types that are adjacent in the ranking, check that the difference between the two polynomials that give the frequencies of those hands never again changes sign.
 
\subsection{Frequency versus showdown frequency}
The Gadbois paradox stems from the fact that in poker with wild cards, not only do hand frequencies affect hand rankings, but hand rankings also affect hand frequencies. In the absence of wild cards, hand rankings do not affect their inherent frequency as seven-card hands. Instead, hand rankings affect player choices which affects each hand's \emph{showdown frequency} as a five-card hand.

We define \emph{showdown frequency} in terms of a preexisting hand ranking. The showdown frequency is the observed frequency of a hand type being declared at showdown. The player chooses which hand type to declare based on the cards available to him and the preexisting hand ranking. Of course, the player can choose not only which five of seven cards to use but will also declare the hand type of his choice based on the ranking, even if the five-card hand has multiple types present. In mathematical notation, we have that the showdown frequency, with respect to a preexisting rank $rk$, is given by:

\[
showdown\text{-}freq_r(h;rk)=|\{s\in S_r|h\in s\text{ and }\forall h'>_{rk}h\text{, }h'\not\in s\}|
\]

We determine showdown-frequencies in the setting where $rk$ is the canonical \emph{frequency ranking} previously. We count a set of seven cards as contributing to hand type $h$ if there exists a subset of five cards with hand type $h$ and there does not exist any subset of five cards with a hand type of higher rank than $h$ according to the frequency ranking. The results for short-deck and standard poker are below, presented in order of increasing showdown frequency:

\begin{table}[h!]
\centering
\begin{tabular}{|lr|lr|}
\hline
\multicolumn{2}{|c|}{\begin{tabular}[c]{@{}c@{}}Short deck\\ $r=9$\end{tabular}} & \multicolumn{2}{c|}{\begin{tabular}[c]{@{}c@{}}Standard deck\\ $r=13$\end{tabular}} \\ \hline
\multicolumn{1}{|l|}{SF}                       & 10,560                          & \multicolumn{1}{l|}{SF}                         & 41,584                            \\ \hline
\multicolumn{1}{|l|}{4X}                       & 44,640                          & \multicolumn{1}{l|}{4X}                         & 224,848                           \\ \hline
\multicolumn{1}{|l|}{FL}                       & 175,560                         & \multicolumn{1}{l|}{FH}                         & 3,473,184                         \\ \hline
\multicolumn{1}{|l|}{HC}                       & 233,100                         & \multicolumn{1}{l|}{FL}                         & 4,047,644                         \\ \hline
\multicolumn{1}{|l|}{3X}                       & 607,200                         & \multicolumn{1}{l|}{ST}                         & 6,180,020                         \\ \hline
\multicolumn{1}{|l|}{FH}                       & 633,024                         & \multicolumn{1}{l|}{3X}                         & 6,461,620                         \\ \hline
\multicolumn{1}{|l|}{ST}                       & 1,169,940                       & \multicolumn{1}{l|}{HC}                         & 23,294,460                        \\ \hline
\multicolumn{1}{|l|}{1P}                       & 2,316,600                       & \multicolumn{1}{l|}{2P}                         & 31,433,400                        \\ \hline
\multicolumn{1}{|l|}{2P}                       & 3,157,056                       & \multicolumn{1}{l|}{1P}                         & 58,627,800                        \\ \hline
\end{tabular}
\caption{Showdown frequencies of hands for short deck and standard poker.}
\end{table}

We may go on to construct a \emph{showdown ranking} of hands based on the likelihood of actually seeing a given hand type declared at showdown. These rankings are given by reading the order that hands appear in Table 1; for example when $r=9$, the showdown ranking gives the lowest rank to 2P. Our wild-card-free version of the Gadbois paradox is when the showdown ranking (with respect to the frequency ranking) is not the same as the frequency ranking itself.\footnote{Note that one may recursively consider the showdown ranking with respect to a ranking that is itself a showdown ranking. For some values of $r$, this process, when initiated with the frequency ranking, creates a forever alternating pair of two rankings. In other cases, the ranking converges. The asymptotics provided in Table 2 can be used to demonstrate that as $r\rightarrow\infty$, the frequency ranking and all recursively defined showdown rankings are equal.}

A few observations about the showdown ranking:

\begin{itemize}
    \item In short deck poker, the showdown ranking is a radical departure from the frequency ranking. Counterintuitively, 1P ranks above 2P, and 3X is rarer than a FH even though a 2P and a FH contain a 1P and a 3X, respectively. In particular, when one conditions on holding a hand with three-of-a-kind, one will expect to improve more often than not.

    \item In a standard deck, if we set aside the likelihood of showing down with merely a high-card hand, then the showdown ranking is the same as the frequency ranking. Moreover, one can check that the standard deck is the deck with the fewest ranks with this property that frequency ranking and showdown ranking agree other than in the case of high-card hands.

    \item The smallest $r$ such that frequency ranking and showdown ranking are equal (including consideration of high-card hands) is $r=23$.
\end{itemize}

\section{Conclusion}
We have provided frequency-based hand rankings for poker played with a card deck with $r$ ranks and four suits, for all $r\geq 5$. Perhaps surprisingly, the hand rankings do not stabilise until $r$ is relatively large. We also introduced the notion of showdown ranking, connecting it to a version of the Gadbois paradox that exists even without the introduction of wild cards.  

\begin{table}[h!]
\renewcommand{\arraystretch}{1.2}
\begin{tabular}{|l|l|r|}
\hline
Hand type & Frequency ranking & Asymptotics \\ 
\hline
SF($r\geq 6$) & 
\begin{tabular}[c]{@{}l@{}}
$\binom{4}{1}\binom{2}{1}\binom{4r-6}{2}+\pmb{\binom{4}{1}\binom{r-5}{1}\binom{4r-7}{2}}$ \\[3pt]
$+\binom{4}{1}\binom{2}{1}\binom{4r-7}{1}+\binom{4}{1}\binom{r-6}{1}\binom{4r-8}{1}$ \\[3pt]
$+\binom{4}{1}\binom{r-5}{1}$
\end{tabular} & $32r^3$ \\ 
\hline
ST($r\geq 8$) & 
\begin{tabular}[c]{@{}l@{}}
$\pmb{\binom{r-5}{1}\cdot 4^5\cdot\binom{4r-28}{2}}+\binom{2}{1}\cdot 4^5\cdot\binom{4r-24}{2}$ \\[3pt]
$+\binom{r-5}{1}\binom{5}{1}\binom{4}{2}4^4\binom{4r-28}{1}+\binom{2}{1}\binom{5}{1}\binom{4}{2}4^4\binom{4r-24}{1}$ \\[3pt]
$+\binom{r-5}{1}\binom{5}{2}\binom{4}{2}^24^3+\binom{2}{1}\binom{5}{2}\binom{4}{2}^24^3$ \\[3pt]
$+\binom{r-3}{1}\binom{5}{1}\binom{4}{3}4^4$ \\[3pt]
$+\binom{r-6}{1}4^6\binom{4r-32}{1}+\binom{2}{1}4^6\binom{4r-28}{1}$ \\[3pt]
$+\binom{r-6}{1}\binom{6}{1}\binom{4}{2}4^5+\binom{2}{1}\binom{6}{1}\binom{4}{2}4^5$ \\[3pt]
$+\binom{r-7}{1}4^7+\binom{2}{1}4^7$
\end{tabular} & $8192r^3$ \\ 
\hline
4X & $\pmb{\binom{r}{1}\binom{4r-4}{3}}$ & $\frac{32}{3}r^4$ \\ 
\hline
FH & 
\begin{tabular}[c]{@{}l@{}}
$\binom{r}{1}\binom{r-1}{1}\binom{4}{2}\binom{r-2}{1}\binom{4}{1}$ \\[3pt]
$+\frac{1}{2}\binom{r}{1}\binom{4}{3}\binom{r-1}{1}\binom{4}{3}\binom{r-2}{1}\binom{4}{1}$ \\[3pt]
$+\frac{1}{2}\binom{r}{1}\binom{4}{3}\binom{r-1}{1}\binom{4}{2}\binom{r-2}{1}\binom{4}{2}$ \\[3pt]
$+\pmb{\frac{1}{2}\binom{r}{1}\binom{4}{3}\binom{r-1}{1}\binom{4}{2}\binom{r-2}{1}\binom{4}{1}\binom{r-3}{1}\binom{4}{1}}$ \\[3pt]
$+\binom{r}{1}\binom{r-1}{1}\binom{4}{3}$
\end{tabular} & $192r^4$ \\ 
\hline
3X & 
\begin{tabular}[c]{@{}l@{}}
$2\binom{r}{2}\binom{4}{3}+\binom{r}{3}\binom{3}{1}\binom{4}{3}\binom{4}{3}\binom{4}{1}$ \\[3pt]
$+\binom{r}{4}\binom{4}{2}\binom{2}{1}\binom{4}{3}\binom{4}{2}\binom{4}{1}\binom{4}{1}+\binom{r}{3}\binom{3}{1}\binom{4}{3}\binom{4}{2}\binom{4}{2}$ \\[3pt]
$+\pmb{\binom{r}{5}\binom{5}{1}4^5}+6\binom{r}{3}\binom{4}{2}\binom{4}{1}+\binom{r}{4}\binom{4}{1}4^3$
\end{tabular} & $\frac{128}{3}r^5$ \\ 
\hline
2P & 
\begin{tabular}[c]{@{}l@{}}
$2\binom{r}{2}\binom{4}{3}+\binom{r}{3}\binom{3}{1}\binom{4}{3}\binom{4}{3}\binom{4}{1}$ \\[3pt]
$+\binom{r}{4}\binom{4}{2}\binom{2}{1}\binom{4}{3}\binom{4}{2}\binom{4}{1}\binom{4}{1}+\binom{r}{3}\binom{3}{1}\binom{4}{3}\binom{4}{2}\binom{4}{2}$ \\[3pt]
$+6\binom{r}{3}\binom{4}{2}\binom{4}{1}+\binom{r}{4}\binom{4}{1}4^3$ \\[3pt]
$+\binom{r}{4}\binom{4}{1}\binom{4}{2}^3\binom{4}{1}+\pmb{\binom{r}{5}\binom{5}{2}\binom{4}{2}^2\binom{4}{1}^3}$
\end{tabular} & $192r^5$ \\ 
\hline
1P & $\pmb{\binom{4r}{7}-4^7\binom{r}{7}}$ & $\frac{256}{5}r^6$ \\ 
\hline
FL & $\pmb{\binom{4}{1}\binom{r}{5}\binom{3r}{2}+\binom{4}{1}\binom{r}{6}\binom{3r}{1}+\binom{4}{1}\binom{r}{7}}$ & $\frac{211}{1260}r^7$ \\ 
\hline
NO & $\pmb{|\{\text{no pair}\}|}-|\text{ST}|+|\text{ST} \cap \text{1P}|-\pmb{|\text{FL}|}+|\text{FL} \cap \text{1P}|+|\text{ST}\cap\text{FL} - \text{1P}|$ & $\frac{37}{12}r^7$ \\ 
\hline
\end{tabular}
\caption{Frequency rankings of poker hands and asymptotics. Bold terms are asymptotically relevant terms.}
\end{table}
\label{table:asymptotics}

\bibliography{pokerReferences}

\end{document}